\documentclass{amsart}
\usepackage{color}

\usepackage{amsmath} 
\usepackage{amssymb}
\usepackage{mathrsfs}

\newtheorem{theorem}{Theorem}[section]

\theoremstyle{definition}
\newtheorem{definition}[theorem]{Definition}

\newtheorem{fact}[theorem]{Fact}

\theoremstyle{remark}
\newtheorem{remark}[theorem]{Remark}

\newtheorem{notation}[theorem]{Notation}

\newcommand{\dn}{{\rm dn}}
\newcommand{\md}{{\rm md}}

\newcommand{\Max}{{\rm Max}}

\newcommand{\lt}{{\rm lt}}
\newcommand{\nacc}{{\rm nacc}}
\newcommand{\acc}{{\rm acc}}

\newcommand{\up}{{\rm up}}

\newcommand{\Rang}{{\rm Rang}}
\newcommand{\lqq}{{{`}{`}}}

\newcommand{\rest}{{\restriction}}
\newcommand{\dom}{{\rm dom}}

\newcommand{\wilog}{{\rm without loss of generality}}

\newcommand{\then}{{\underline{then}}}

\newcommand{\when}{{\underline{when}}}

\newcommand{\mn}{{\medskip\noindent}}
\newcommand{\sn}{{\smallskip\noindent}}

\newcommand{\bfc}{{\mathbf c}}

\newcommand{\bfT}{{\mathbf T}}

\newcommand{\cP}{{\mathscr P}}

\newcommand{\varp}{{\varepsilon}}

\newcommand{\cU}{{\mathscr U}}

\newcommand{\extheta}{{\partial }}
\newcommand{\W}{{{\mathscr W } }}

\newcommand{\cf}{{\rm cf}}

\newcount\skewfactor
\def\mathunderaccent#1#2 {\let\theaccent#1\skewfactor#2
\mathpalette\putaccentunder}
\def\putaccentunder#1#2{\oalign{$#1#2$\crcr\hidewidth
\vbox to.2ex{\hbox{$#1\skew\skewfactor\theaccent{}$}\vss}\hidewidth}}

\newenvironment{PROOF}[2][\proofname.]
   {\begin{proof}[#1]\renewcommand{\qedsymbol}{\textsquare$_{\rm #2}$}}
   {\end{proof}}

\usepackage[hidelinks]{hyperref}

\begin{document}
\makeatletter\def\shfiuwefootnote{\gdef\@thefnmark{}\@footnotetext}\makeatother\shfiuwefootnote{Version 2021-05-06\_5. See \url{https://shelah.logic.at/papers/1163/} for possible updates.}

\title {Colouring of successor of regular again \\
Sh1163}
\author {Saharon Shelah}
\address{Einstein Institute of Mathematics\\
Edmond J. Safra Campus, Givat Ram\\
The Hebrew University of Jerusalem\\
Jerusalem, 91904, Israel\\
 and \\
 Department of Mathematics\\
 Hill Center - Busch Campus \\ 
 Rutgers, The State University of New Jersey \\
 110 Frelinghuysen Road \\
 Piscataway, NJ 08854-8019 USA}
\email{shelah@math.huji.ac.il}
\urladdr{http://shelah.logic.at}
\thanks{Paper number 1163.
First version on May 17, 2019.  References like
\cite[Th.2.2=Le8]{Sh:1127} means the label of Th.2.2 is e8.
The reader should note that the
version in my website is usually more updated than the one in the
mathematical archive}







\date{2021-05-06d}

\begin{abstract}        




We get a  version of the colouring property $\Pr_1$  proving
$\Pr_1(\lambda,\lambda,\lambda,\partial)$ 
always when $\lambda
= \partial^+,\partial   $ 
are regular cardinals and some stationary subset of $ \lambda $ 
consisting of ordinals of cofinality $ < \partial $ do not
reflect in any ordinal $ < \lambda $.  


 
 
 






\end{abstract}

\maketitle
\numberwithin{equation}{section}
\setcounter{section}{-1}
\newpage

\section {Introduction}

 and 


We prove a strong colouring theorem on successor of regular
uncountable cardinals, so called $ \Pr_1$.  

 On the history of $ \Pr_1 $
see  \cite[Ch.III,\S4]{Sh:g} 
and later \cite{Sh:572}, and then independently 
Rinot \cite{Ri14} and \cite{Sh:1027}.

Rinot \cite[Main result]{Ri14} proved that $ \Pr_1(\lambda, \lambda, \lambda, \theta )$  
        when those are regular cardinals;
   $ \lambda = \theta ^{++} $    or just  
        $\theta ^+ < \lambda $ 
        and $ \lambda $ is a successor of 
        regular or just it has a 
        non-reflecting stationary subset
        of $ \lambda $ consisting of ordinals of cofinality 
        at least $ \theta $.
    In \cite{Sh:1027}, 
        we have $ \Pr_1( \lambda ,   
            \lambda, \lambda , (\theta _0, \theta ))$   
            where   
            $ \theta _0$  
            is regular $ < \theta = \cf( \theta ) , \theta ^+ < \lambda $ and $ \lambda $ is a successor of regular.
             
              Earlier  \cite[4.2,
            page 27]{Sh:572}  
      prove that 
            $ \Pr_1( \lambda , \lambda , \lambda , \theta ) $ when in addition
            $ \lambda = \theta ^{++}$.
    
    Much earlier \cite[Ch.III, \S4]{Sh:g}    
        had treated   
those problems  
            in a general but probably not 
            so transparent  
            way, first 4.1 gives a set of various hypothesis
            (each with some parameters) 

The result    here  
is incomparable with the ones  
in \cite{Ri14}, \cite{Sh:1027}, \cite{Sh:572}  
the assumption on the stationary set is stronger but
the the arity - the last parameter, $ \theta $  
is bigger.

The connection between purely combinatorial theorems and topological
constructions is known for many years. 
Several
results in general topology were proved using the property
$\Pr_1(\lambda,\mu,\sigma,\theta)$, 
see recently \cite{Sh:1025}, then \cite[\S1]{Sh:1027}.  

\noindent
Recall: 
\begin{definition}
\label{e18}
1) Assume $\lambda \ge \mu \ge \sigma + \theta_0 + \theta_1,\bar\theta =
(\theta_0,\theta_1)$, see \ref{e28}(1).  
Assume further that $ \theta _0, \theta _1  \ge {\aleph_0} $ but $ \sigma $
 may be finite

Let $\Pr_1(\lambda,\mu,\sigma,\bar\theta)$ mean
that there is $\mathbf c:[\lambda]^2 \rightarrow \sigma$ witnessing it,
which means:
\mn
\begin{enumerate}
\item[$(*)_{\bfc}$]  if (a) then (b), where:
\sn
\begin{enumerate}
\item[(a)]  for $\iota=0,1,\mathbf i_\iota < \theta_\iota$ and
  $\bar\zeta^\iota = \langle \zeta^\iota_{\alpha,i}:\alpha < \mu,i <
  \mathbf i_\iota\rangle$ are sequences of ordinals of $\lambda$ without
  repetitions,
  and 
  $\Rang(\bar\zeta^0) ,  $ 
   $\Rang(\bar\zeta^1 )  $  
  are disjoint and $\gamma < \sigma$
\sn
\item[(b)]  there are $\alpha_0 < \alpha_1 < \mu$ such that $\forall
  i_0 < \mathbf i_0,\forall i_1 < \mathbf i_1,\mathbf c
  \{\zeta^0_{\alpha_0,i_0},\zeta^1_{\alpha_1,i_1}\} = \gamma$ 
  and $  \zeta^0_{\alpha_0,i_0}  < \zeta^1_{\alpha_1,i_1} $.
\end{enumerate}
\end{enumerate}
\mn
2) Above if $\theta_0 = \theta = \theta_1$ then we may write
$\Pr_1(\lambda,\mu,\sigma,\theta)$. 
\end{definition}


In this paper  we prove e.g. that  
    if some stationary $ S \subseteq \{\delta < { \aleph_2 } 
    :  \cf( \delta ) < {\aleph_1} \} $ do not reflect then  
$\Pr_1(\aleph_2,\aleph_2,\aleph_2,\aleph_1)$ holds,
which means that countable infinite sequences can be taken in both
$\lqq$sides". 
Actually, the theorem says that, in particular, 
$\Pr_1(\lambda,\lambda,\lambda,\extheta)$ holds whenever $\extheta =
\cf(\extheta)$ and $\lambda = \extheta^{+}$ 
    and there is 
a non-reflecting stationay subset of $ \lambda $.  
  We intend to say more on other $ \lambda $-s in \cite{Sh:F2047}.

We thank the referee for many good suggestions.

\begin{definition}
\label{e24}
1)
A filter $D$ on a set $I$ is
uniform \when \, for every subset $ A $ of $ I $ of 
cardinality   
$ < |I|$, the set $ I \setminus A \in D$; all our 
filters will be uniform

\noindent 
2)
A filter $D$ on a set $I$ is
weakly $\theta$-saturated \when \,
$ \theta \ge |I|  $  and 
there
is no partition of $I$ to $\theta$ sets from $D^+$, 

\noindent 
3)
We say the filter $ D$ on a set $ I $ is   $\theta$-saturated \when \, 
the Boolean algebra  
$\cP(I)/D$ satisfies the $\theta$-c.c.

\end{definition}

\noindent
\begin{fact}
\label{e25}
1)
If $D$ is a $\theta$-complete filter on $\lambda$ and is not
$\theta$-saturated \then  \, 
it is not weakly $\theta$-saturated; 
so those properties are equivalent.  

\noindent 
2)
If $\theta = \sigma^+$ and $D$ is a $\theta$-complete filter on
$\theta$, \then \, $D$ is not weakly $\theta$-saturated.

\noindent
3) If $ n \ge 1$  and $ \lambda = \sigma ^{+ n }$ and $ D $  
is a (uniform) $ \sigma ^+ $-complete filter on 
$ \lambda $ \then \, $ D $ is not weakly 
$ \sigma ^{+n} $-saturated  

\end{fact}

\begin{PROOF}{\ref{e25}}
1)
Obvious and well  known

\noindent
2) By \cite{So2},  

\noindent  
3) Let $ \mu $ be the minimal cardinal such that 
$ D $ is not $ \mu ^+$-complete, so clearly 
$ \mu \in [ \sigma ^+, \lambda ]$  hence
$ \mu $ is a successor cardinal. 
So there is a function $ f $ from $ \lambda $ into
$ \mu $ such that for every subset $ A$  of $ \mu $ 
of cardinality $ < \mu $, $ f^{-1}(A) = \emptyset \mod D $.
Let $ E $ be the family of subsets $ A $ of $ \mu $ 
such that $f^{-1}(A) \in D $. Clearly $ E $ is a (uniform)
$ \mu $-complete filter on $ \mu $ hence  
by part (2) is not 
weakly $ \mu $-saturated, let $ \langle A_\varepsilon :
\varepsilon < \mu \rangle $ be a partition of $ \mu $
to sets 
from $ E^+ $. Now $\langle f^{-1}(A_\varepsilon ): 
\varepsilon < \mu \rangle $ witnesses the desired conclusion.



\end{PROOF}

\begin{notation}
\label{e28}
1) We denote infinite cardinals by 
$ \lambda, \mu,\kappa,  \theta, \partial $
while $ \sigma $ denotes a finite or infinite cardinal. 
We denote ordinals by $ \alpha , \beta, \gamma , \varepsilon ,
\zeta, \xi $.
Natural numbers are denoted by $ k,{\ell}, m , n $ 
and $ \iota \in \{ 0,1,2 \} $

\noindent 
1A) Let $ D $ denote a filter on an infinite set 
$ \dom(D)$ 

\noindent 
2) For a set $ A $ of ordinals let
$ \nacc (A) = \{\alpha \in A:   
\alpha > \sup(A \cap \alpha ) \} $
and $ \acc (A) = A \setminus \nacc(A)$.  
\noindent  
For regular
    cardinals 
$\lambda > \kappa$ let $S^\lambda_\kappa = \{\delta <
\lambda:\cf(\delta) = \kappa\}$.
\end{notation}
\newpage

\section {A colouring theorem}

\noindent
Our aim is to prove
\begin{theorem}  
\label{f4}
$\Pr_1(\lambda,\lambda,\extheta,\extheta)$  
and 
moreover  
    $ \Pr_1(\lambda,\lambda,\lambda,\extheta)$ holds 
provided that:
\mn
\begin{enumerate}
\item[(a)]  $\lambda = \partial^+ $
\sn
\item[(b)] $ \partial = \cf(\partial)   > \aleph_0$
  
\sn
\item[(c)]    
  $ {\mathscr W } $ is a stationary subset 
  of $ \lambda $
consisting of ordinals of cofinality $ < \partial $  reflecting
in no ordinal  $ < \lambda $   

\end{enumerate} 
\end{theorem}

\begin{remark}  
\label{f7}  
1)
The case of $\extheta$ colours, i.e. proving only
$\Pr_1(\lambda,\lambda,\extheta,\extheta)$  
is easier so we prove it first.





 

\noindent 
2)  
Can we weaken clause (c)  
of \ref{f4} replacing 
$\lqq$reflecting in no ordinal $ < \lambda $ by
$ \lqq $reflecting in no ordinal of cofinality $  \partial $?  

The answer seem yes provided that we add:  

\begin{enumerate} 
\item[$(
    \alpha ) $]  there is a sequence $ \langle e_ \alpha: 
    \alpha \notin \W\rangle $ such that 
    ($ \W $ is as above and)
    $ e_ \alpha $
    is a club of $ \alpha $ of order type $ < \partial $ 
    and  for  $ \alpha \in   
        e_ \beta $ from $ \W$  
    we have $ e_\alpha = \alpha \cap e_ \beta $
\item[$(\beta )$]  
there is no $\extheta$-complete not 
$\extheta^+$-complete uniform weakly
$\extheta$-saturated filter on $\lambda$.
\end{enumerate} 
\end{remark}

\begin{PROOF}{\ref{f4}}  
\medskip

\noindent
\underline{Stage A}:  We begin   
as in earlier proofs (e.g.  \cite{Sh:1027}).
We let
$ (\kappa_1,\kappa_2) = 
(\extheta ,\lambda)$.
Let $S \subseteq  
      S^\lambda_\partial$ 
be stationary and $h:\lambda
\rightarrow \lambda$ be such that $\alpha < \lambda \Rightarrow
h(\alpha) < 1 + \alpha,h \rest (\lambda \backslash S)$ is
constantly zero and $S^*_\gamma := \{\delta \in S:h(\delta) = \gamma\}$
is a stationary subset of $\lambda$ for every $\gamma < \lambda$.  Let
$F_\iota:\lambda \rightarrow \kappa_\iota$ for $\iota = 
1,2$ be such
that for every 
$(  \varepsilon_1,\varepsilon_2) \in
 (\kappa_1 \times \kappa_2)$ the set
$W_  
{\varepsilon_1,\varepsilon_2}(\beta) = \{\gamma \in
S^*_\beta:F_\iota(\gamma) = \varepsilon_\iota$ for 
$\iota = 1,2 \} $   
is a
stationary subset of $\lambda$ for every $\beta < \lambda$.

For $\iota = 
1,2$ and 
$\rho \in {}^{\omega >}\lambda$ let $F_\iota(\rho) = \langle
F_\iota(\rho(\ell)):\ell < \ell g(\rho)\rangle$.

\begin{enumerate} 
\item[$ \odot_0$]  \wilog \,  
    if $ \delta \in \W $  then  
    $ \delta $ is divisible
    by $ \partial $.
\end{enumerate} 

Let $\bar e = \langle e_\alpha:\alpha < \lambda\rangle$ be such that:
\mn
\begin{enumerate}
\item[$\odot_1$]
\begin{enumerate}
\item[(a)]   if $\alpha = 0$ then $e_\alpha = \emptyset$
\sn
\item[(b)]   if $\alpha = \beta +1$ then $e_\alpha = \{\beta\}$
\sn
\item[(c)]   if $\alpha$ is a limit ordinal then $e_\alpha$ is a club of
  $\alpha$ of order type $\cf(\alpha)$ disjoint
to $S^\lambda_\partial$ hence to $S$.
\item[(d)]  if $ \alpha $   
    is a limit ordinal 
   then $ e_\alpha $ 
    is disjoint to $ \W$.  
\end{enumerate}
\end{enumerate}
\mn
In other cases 
(not here)  
instead $ h $   we use a sequence 
$ \langle  h_\alpha : \alpha < \lambda \rangle $ of functions,  
$h_ \alpha : e_ \alpha \rightarrow \extheta   $
and use  e.g 
$ \langle h_ {\gamma _ {\ell}(\beta, \alpha )}
  ( \gamma _{{\ell} + 1}(\beta, \alpha ) ) : 
  {\ell} < k( \beta, \alpha )\rangle $  
  and $ \rho _{ h }$, but   
  this is 
  not necessary here.
  

Now (using $\bar e$) for   
$ \alpha < \beta < \lambda$, let

\[
\gamma(\beta,\alpha) := \min\{\gamma \in e_{\beta}:\gamma \ge
\alpha\}.
\]

\mn
Let us define $\gamma_{\ell}(\beta,\alpha)$:

\[
\gamma_{0}(\beta,\alpha) = \beta,
\]

\[
\gamma_{\ell +1}(\beta,\alpha) = \gamma(\gamma_{\ell}(\beta,\alpha),
\alpha) \text{ (if  well 
defined)}.
\]

\mn
If 
$ \alpha < \beta < \lambda$, let $k(\beta,\alpha)$ be the maximal $k < 
\omega$ such that $\gamma_{k}(\beta,\alpha)$ is defined (equivalently
is equal to $\alpha$) and let
$\rho_{\beta,\alpha} = \rho(\beta,\alpha)$  be the sequence  

\[
\langle \gamma _{0}(\beta,\alpha),\gamma_{1}(\beta,\alpha),\ldots,
\gamma_{k(\beta,\alpha)-1}(\beta,\alpha)\rangle.
\]

\mn
Let $\gamma_{\ell t}(\beta,\alpha) =
\gamma_{k(\beta,\alpha)-1}(\beta,\alpha)$ where $\ell t$ stands for last.   

Let

\[
\rho _h = \langle h(\gamma _{\ell} (\beta, \alpha )): 
   {\ell} < k( \beta, \alpha )\rangle 
\]

\mn
and we let $\rho(\alpha,\alpha)$ and  $\rho_{h}(\alpha,\alpha)$ 
be the empty sequences.  
\noindent
Now clearly:
\mn 
\begin{enumerate}
\item[$\odot_2$]   if 
$ \alpha < \beta < \lambda$ then $\alpha \le
\gamma(\beta,\alpha) < \beta$
\end{enumerate}
\mn
hence
\mn 
\begin{enumerate}
\item[$\odot_3$]   if 
$ \alpha < \beta < \lambda,0 < \ell < \omega$, and
$\gamma_{\ell}(\beta,\alpha)$ is well defined, then

\[
\alpha \le \gamma_{\ell}(\beta,\alpha) < \beta 
\]
\end{enumerate}
\mn
and
\mn 
\begin{enumerate}
\item[$\odot_4$]  if 
$ \alpha < \beta < \lambda$, then $k(\beta,\alpha)$ is
well defined and letting $\gamma_{\ell} :=
\gamma_{\ell}(\beta,\alpha)$ for $\ell \le k(\beta,\alpha)$ we have
 
\[
\alpha = \gamma_{k(\beta,\alpha)} < \gamma_{\ell t}(\beta,\alpha) = 
\gamma_{k(\beta,\alpha)-1} < \cdot \cdot \cdot < \gamma_{1} < \gamma_{0} = 
\beta 
\]

\[
\text{ and } \alpha \in e_{\gamma_{\ell t}(\beta,\alpha)}
\]

\mn
i.e. $\rho(\beta,\alpha)$ is a (strictly) decreasing finite 
sequence of ordinals, starting with $\beta$, ending with 
$\gamma_{\ell t}(\beta,\alpha)$ of length $k(\beta,\alpha)$.
\end{enumerate}
\mn
Note that if $\alpha \in S,\alpha < \beta$ then
$\gamma_{\ell t}(\beta,\alpha) = \alpha +1$.

Also
\mn 
\begin{enumerate}
\item[$\odot_5$]  if $\delta$ is a limit ordinal
 and $\delta < \beta < \lambda$,
\then \, for some $\alpha_{0} < \delta$ we have: $\alpha_{0} \le \alpha <
\delta$ \underline{implies}:
\sn
\begin{enumerate}
\item[$(i)$]  for $\ell < k(\beta,\delta)$ we have 
$\gamma_{\ell}(\beta,\delta) = \gamma_{\ell}(\beta,\alpha)$   
\sn
\item[$(ii)$]  $\delta \in \nacc(e_{\gamma_{\ell t}(\beta,\delta)})
\Leftrightarrow \delta = \gamma_{k(\beta,\delta)}(\beta,\delta) =
\gamma_{k(\beta,\delta)}(\beta,\alpha) \Leftrightarrow 
\neg[\gamma_{k(\beta,\delta)}(\beta,\delta) = \delta > 
\gamma_{k(\beta,\delta)}(\beta,\alpha)]$   
\sn
\item[$(iii)$]  $\rho(\beta,\delta) \trianglelefteq
  \rho(\beta,\alpha)$; i.e. is an initial segment
\sn
\item[$(iv)$]  $\delta \in \nacc(e_{\gamma_{\ell t}(\beta,\delta)})$
  (here always holds if $\delta \in S$) implies:
\begin{itemize}
\item  $\rho(\beta,\delta) \char 94 \langle \delta \rangle \trianglelefteq 
\rho(\beta,\alpha)$ hence 
\sn
\item  $\rho_{h}(\beta,\delta) \char 94 
\langle h 
  (\beta,\delta)
(\delta )\rangle \trianglelefteq \rho_{h}(\beta,\alpha)$.
\end{itemize}
\sn
\item[$(v)$]  if $\cf(\delta) = \partial$ 
    or $ \delta \in {\mathscr W } $  
\then \, we have 
$\gamma_{\ell t}(\beta,\delta) = \delta +1$ so $\delta 
    + 1 \in  2021-05-06   
\nacc(e_{\gamma_{\lt}(\beta,\delta)})$ 
\sn
\item[$(vi)$]  if $\cf(\delta) = \partial$ 
        or $ \delta \in {\mathscr W } $ 
and 
    $\delta \in e_\gamma$
     or $ \delta \in {\mathscr W } $,  
\then \, necessarily $\gamma = \delta +1$.
\end{enumerate}
\end{enumerate}
\mn
Why?  Just let

\[
\alpha_{0} = \Max\{\sup(e_{\gamma_{\ell}(\beta,\delta)} \cap \delta) +
1: \ell < k(\beta,\delta) \text{ and } \delta \notin
\acc(e_{\gamma_{\ell}(\beta,\delta)})\}.
\]

\mn
Notice that if $\ell < k(\beta,\delta)-1$ then $\delta \notin
\acc(e_{\gamma_\ell(\beta,\delta)})$ follows.

Note that the outer maximum (in the choice of $\alpha_0$) 
is well defined as it is over a finite non-empty set of 
ordinals.  The inner $\sup$ is on the empty set (in which case we get
zero) or is the maximum
 (which is well defined) as $e_{\gamma_{\ell}(\beta,\delta)}$ is a 
closed subset of $\gamma_{\ell}(\beta,\delta),\delta < \gamma_{\ell}
(\beta,\delta)$ and $\delta \notin
\acc(e_{\gamma_{\ell}(\beta,\delta)})$ - as this is required.  For
clauses (v), (vi) recall
    $\delta \in S^\lambda_\partial  \cup {\mathscr W } $  
    and $e_\gamma \cap
S^\lambda_\partial = \emptyset$
    and $ e_ \gamma \cap {\mathscr W } = \emptyset $  
when 
   $\gamma$ is a limit ordinal and
$e_\gamma = \{\gamma - 1    
  \}$ when $\gamma$ is a successor ordinal.
\mn
\begin{enumerate}
\item[$\odot_6$]  
\begin{enumerate}
\item[(a)]  if 
  $ \alpha < \beta < \lambda,
\ell < k(\beta,\alpha),\gamma =
\gamma_{\ell}(\beta,\alpha)$ \then \, $\rho(\beta,\alpha) =
\rho(\beta,\gamma) \char 94 \rho(\gamma,\alpha)$ and 
$\rho_{h}(\beta,\alpha) = 
\rho_{h}(\beta,\gamma) \char 94 \rho_{h}(\gamma,\alpha)$
\sn
\item[(b)]  if 
$ \alpha_0 < \ldots < \alpha_k$ and
  $\rho(\alpha_k,\alpha_0) = \rho(\alpha_k,\alpha_{k-1}) \char 94
  \ldots \char 94 \rho(\alpha_1,\alpha_0)$ \then \, this holds for any
  sub-sequence of $\langle \alpha_0,\dotsc,\alpha_k\rangle$.
\end{enumerate}
\sn
\item[$\odot_7$]  let $F'_\iota $ 
    be
  $F_\iota \circ h$ for $\iota = 1,2$; so   
  $ F'_1 $ is a function from $ \lambda $ into $ \extheta $ 
  and $ F'_2 $ is a function from $ \lambda $ into $ \lambda $.
    
\end{enumerate}
\medskip

\noindent
\underline{Stage B}:

Let
\mn
\begin{enumerate}
\item[$\boxplus_2$]  $\mathbf T = \{\bar t:\bar t = \langle
t_\alpha:\alpha < \lambda\rangle$ satisfies $t_\alpha \in [\lambda]^{<
  \extheta}$ and\footnote{
  if instead we demand $ \alpha \not= \beta < \lambda \Rightarrow 
  t_ \alpha \cap t_\beta = \emptyset $  then we shall get the 
  same filter $ D $.
  } 
  $t_\alpha \subseteq \lambda \backslash \alpha\}$.
\sn
\item[$\boxplus_3$]  for $\varp < \extheta$ 
and $\bar t \in \mathbf T$ let $A_{\bar t,\varepsilon}$ be 
the set of $\gamma < \lambda$ such that for some 
$(\alpha_0,\alpha_1)$ we have: 
\sn
\begin{enumerate}
\item[$(a)$]  $\alpha_0 < \alpha_1 < \lambda$ and\footnote{If we
choose to add here ``$t_{\alpha_0} \subseteq \alpha_1$", \then \, 
we would 
a problem
in proving clause $\boxplus_5(b)$.} $(\zeta,\xi) \in t_{\alpha_0}
\times t_{\alpha_1} \Rightarrow \zeta < \xi$
\sn
\item[$(b)$]  for every $(\zeta,\xi) \in t_{\alpha_0} \times t_{\alpha_1}$
  for some $\ell$ we have:
\sn
\begin{enumerate}
\item[$(\alpha)$]  $\ell < k(\xi,\zeta)$
\sn
\item[$(\beta)$]  $\gamma_\ell(\xi,\zeta) = \gamma$
\sn
\item[$(\gamma)$]   if $k < k(\xi,\zeta)$ then $F'_1(\gamma)
 \ge F'_1(\gamma_k(\xi,\zeta))$ and $F'_1(\gamma) \ge \varepsilon$
\sn
\item[$(\delta)$]   if $k < \ell$ then $F'_1(\gamma_k(\xi,\zeta)) <
  F'_1(\gamma)$.
\end{enumerate}
\end{enumerate}
\end{enumerate}
\mn
We define:
\mn
\begin{enumerate}
\item[$\boxplus_4$]  $D = \{A \subseteq \lambda:A$ includes
$A_{\bar t,\varp}$ for some $\bar t \in \mathbf T,\varepsilon <
\extheta\}$.
\end{enumerate}
\mn
Now note:
\mn
\begin{enumerate}
\item[$\boxplus_5$]
\begin{enumerate}
\item[(a)]   if $\bar s,\bar t \in \bfT,\varp \le \zeta < 
\extheta$ and $(\forall \alpha < \lambda)(s_\alpha \subseteq t_\alpha)$, 
\then \, $A_{\bar t,\zeta} \subseteq A_{\bar s,\varp}$
\sn
\item[(b)]  if $\bar s \in \bfT,\varp < \extheta,g$ is an increasing
function from $\lambda$ to $\lambda$ and $\bar t = \langle
  t_\alpha:\alpha < \lambda\rangle$ is defined by $t_\alpha = s_{g(\alpha)}$
  \then \, $A_{\bar t,\varp} \subseteq A_{\bar s,\varp}$.
\end{enumerate}
\end{enumerate}
\mn
[Why?  Read the definitions.]
\mn
\begin{enumerate}
\item[$\boxplus_6$]  
\begin{enumerate}
\item[(a)]   the intersection of any $< \extheta$ members of $D$ is
 a member of $D$, equivalently includes the 
set $A_{\bar t,\zeta}$ for some $\bar t \in \mathbf T,\zeta < \extheta$
\sn
\item[(b)]  for every $\beta < \lambda$ for some 
$\bar t \in \mathbf T,A_{\bar t,0} \subseteq [\beta,\lambda)$
\sn
\item[(c)]  if $\bar t \in \bfT$ and $\alpha < \lambda \Rightarrow
t_\alpha \ne \emptyset$ then $\cap\{A_{\bar t,\varp}:\varp < \extheta\} 
= \emptyset$
\sn  
\item[(d)]  $D$ is upward closed.
\sn
\item[(e)]   $ \lambda $ belongs to $ D $  
\end{enumerate}
\end{enumerate}
\mn
[Why?  For clause (a) assume $A_\varepsilon \in D$ for 
$\varepsilon < \varepsilon(*) < \extheta$ then for some $\zeta_\varepsilon <
\extheta$ and $\bar t_\varepsilon 
\in \mathbf T$ we have $A_\varepsilon \supseteq 
A_{\bar t_\varepsilon,\zeta_\varepsilon}$.  Define $t_\alpha =
 \bigcup\{t^\varepsilon_\alpha:\varepsilon < \varepsilon(*)\}$ 
for $\alpha < \lambda$ and 
$\zeta = \sup\{\zeta_\varepsilon:\varepsilon < \varepsilon(*)\}$;
as  
the cardinal $\extheta$ is regular, clearly
$|t_\alpha| \le \sum\limits_{\varepsilon < \varepsilon(*)}
|t^\varp_\alpha| < \extheta$ and 
obviously 
$t_\alpha \subseteq 
[\alpha,\lambda)$ hence $\bar t = \langle t_\alpha:\alpha <
  \lambda\rangle \in \mathbf T$ and similarly $\zeta < \extheta$.  
Easily $A_{\bar t,\zeta} \subseteq A_{\bar t_\varepsilon,\zeta_\varepsilon}$ 
for every $\varepsilon < \varepsilon(*)$, see $\boxplus_5(a)$ so we
are done proving clause (a).  
For clause (b) define $t_\alpha = 
\{\beta + \alpha +1\}$ and 
recalling $ \boxplus_3(b)(\beta )$  
and $ \odot _4$
check that 
$A_{\bar t,0} \subseteq [\beta,\lambda)$.  Also clause (c) obviously
holds because $\gamma \in A_{\bar t,\varp} \Rightarrow F'_1(\gamma)
\ge \varp$ by $\boxplus_3(b)(\gamma)$ 
and $ F'_1 $  
is a function from $ \lambda $ to $ \extheta $ 
and clauses 
(d),(e)  
hold trivially by the
definition.]
\mn
\begin{enumerate}
\item[$\boxplus_7$]  
\begin{enumerate}
\item[(a)]  $\emptyset \notin D$
\sn
\item[(b)]  $D$ is a filter on $\lambda$, equivalently 
$A_{\bar t,\varepsilon} \ne \emptyset$ for every $\bar t,\varp$; 
also $D$ is uniform $\extheta$-complete, 
not $\extheta^+$-complete.
\end{enumerate}
\end{enumerate}
\mn
[Why?  Clause (a) is a major point, proved in Stage C below.  That is,
by $\boxplus_6(a),(d)$ the only missing 
point is to show $A_{\bar t,\zeta} \ne \emptyset$, (in fact, 
$|A_{\bar t,\zeta}| = \lambda$).  For clause (b) by (a) and
$\boxplus_6(a),(d),(e),D$ 
is a $\extheta$-complete filter and the statement that
$D$ is uniform holds by $\boxplus_6(b)$ and not $\extheta^+$-complete
holds by $\boxplus_6(c)$.]

Note also
\mn
\begin{enumerate}
\item[$\boxplus_8$]  $D$ is not weakly $\extheta$-saturated. 
\end{enumerate}
\mn
[Why?  By $\boxplus_7  + \boxplus_6(c)$   
and clause (c)   
    in the assumptions of the theorem.
That is it is known that if $ D $ fail this statement (and has 
the properties listed before) then there is no 
    $ {\mathscr W } $ as  
in clause (c)   
of the theorem.  
That is, considering the forcing notion $ \mathbb{P} =D^+$    
with inverse inclusion;
in $ \mathbf{V} ^{\mathbb{P} }$, the generic set $ \mathbf{G} $
is an ultrafilter on the Boolean algebra 
    $   {\mathscr P} ( \lambda )^ \mathbf{V}$ 
    and let $ \mathbf{j} $ be the canonical embedding from 
    $ \mathbf{V} $ into the Mostowski collapse of 
    $ \mathbf{V} ^\lambda / \mathbf{G} $ 
    (we are using only functions from $ \mathbf{V} $),
    now the contradiction will be clear. 
    If $ \partial $ 
    is a successor cardinal we can use 
    \ref{e25}(2).   

\medskip

\noindent
\underline{Stage C}:

In this stage we accomplish the proof of the missing point in
$\boxplus_7(a)$ from above,
so we shall prove ``$A_{\bar t,\varp}$ is non-empty (in fact, 
has cardinality $\lambda$)" \when \,:
\mn
\begin{enumerate}
\item[$\boxplus$]   
\begin{enumerate}
\item[(a)]  $t_{\alpha} \subseteq \lambda \backslash \alpha$ for
  $\alpha < \lambda$
\sn
\item[(b)]  $|t_\alpha| < \extheta$
\sn
\item[(c)]  $\varp < \extheta$.
\end{enumerate}
\end{enumerate}
\mn
To start we note that:
\mn
\begin{enumerate}
\item[$(*)_1$]  \wilog \, $t_\alpha \ne \emptyset$ and $\alpha 
< \min(t_\alpha)$.
\end{enumerate}
\mn
[Why?  First, recalling $\boxplus_5(a)$ we can replace $\bar t$ by
$\bar t = \langle t_\alpha \cup\{\alpha\}:\alpha < \lambda\}$, so we
may assume that each $t_\alpha$ is not empty.  Second, let 
$\bar t' = \langle t'_\alpha:\alpha <
\lambda\rangle,t'_\alpha = t_{\alpha +1}$, so easily $\bar t'$
satisfies $(*)_1$ and $A_{\bar t',\varp} \subseteq A_{\bar t,\varp}$
by clause $\boxplus_5(b)$.]

Now
\mn
\begin{enumerate}
\item[$(*)_2$]  we can find $\cU^{\dn}_1,\varepsilon^{\dn}$ such that:
\sn
\begin{enumerate}
\item[$(a)$]  $\cU^{\dn}_1 \subseteq    \W  $ 
    is stationary in $\lambda$, see stage A on $ S^*_0$  
\sn
\item[$(b)$]  $\alpha < \delta \in \cU^{\dn}_1 \Rightarrow t_\alpha
  \subseteq \delta$
\sn
\item[$(c)$]  $\varepsilon^{\dn} < \extheta$
\sn
\item[$(d)$]  if $\delta \in \cU^{\dn}_1$ \then \, for arbitrarily large
$\alpha < \delta$ we have $\zeta \in t_\alpha 
\Rightarrow \Rang(F_1(\rho_{h}(\delta,\zeta)))
  \subseteq \varepsilon^{\dn} < \kappa_1  
  = \extheta$.
\end{enumerate}
\end{enumerate}
\mn
[Why?  Clearly $E_0 = \{\delta < \lambda:\delta$ is a limit ordinal such
  that $\alpha < \delta \Rightarrow t_\alpha \subseteq \delta\}$ is a
  club of $\lambda$.
For every $\delta \in   \W  
    \cap E_0$ and $\alpha < \delta$ we can find
$\varepsilon^{\dn}_{\delta,\alpha}$ as in clauses (c),(d) of $(*)_2$
    (because $ | t_ \delta | < \partial $)  
and so recalling that $\cf(\delta)  < \partial $  
it follows that there is $\varepsilon^{\dn}_\delta$
such that $\delta = \sup\{\alpha <
\delta:\varepsilon^{\dn}_{\delta,\alpha} = 
\varepsilon^{\dn}_\delta\}$.  
Then recalling $\lambda = \cf(\lambda) > \extheta $ we can
choose $\varepsilon^{\dn}$ such that the set $\cU^{\dn}_1 = \{\delta
\in \W   
\cap E_0:\varepsilon^{\dn}_\delta = 
\varepsilon^{\dn}\}$ is stationary.  So $(*)_2$ holds indeed.]
\mn
\begin{enumerate}
\item[$(*)_3$]  We can find $\cU^{\up}_1,\alpha^*_1,
\varepsilon^{\up}$ such that:
\sn
\begin{enumerate}
\item[$(a)$]  $\cU^{\up}_1 \subseteq S^*_0$ is stationary
\sn
\item[$(b)$]  $h \rest \cU^{\up}_1$ is constantly 0, actually follows by (a),
  see Stage A
\sn
\item[$(c)$]  $\alpha^*_1 < \lambda$ 
satisfies $\alpha^*_1 < \min(\cU^{\up}_1)$ and
$\varepsilon^{\up} < \extheta$
\sn
\item[$(d)$]  if $\delta \in \cU^{\up}_1$ and $\alpha \in
[\alpha^*_1,\delta)$ and $\beta \in t_\delta$ \then \,:
\sn
\begin{itemize}
\item  $\rho_{\beta,\delta} \char 94 \langle \delta \rangle
  \trianglelefteq \rho_{\beta,\alpha}$
\sn
\item  $\Rang(F_1(\rho_{h}(\beta,\delta))) 
\subseteq \varepsilon^{\up}$.
\end{itemize}
\end{enumerate}
\end{enumerate}
\mn
[Why?  For every $\delta \in S^*_0 \subseteq S$ and $\zeta \in t_\delta$ let
$\alpha_{1,\delta,\zeta} < \delta$ be such that $(\forall \alpha)
(\alpha \in [\alpha_{1,\delta,\zeta},\delta) \Rightarrow
\rho_{\zeta,\delta} \char 94 \langle \delta \rangle 
\trianglelefteq \rho_{\zeta,\alpha})$, it exists by $\odot_5$ of Stage A.

Let
\mn
\begin{itemize}
\item   $\alpha_{1,\delta} = \sup\{\alpha_{1,\delta,\zeta}:\zeta \in
t_\delta\}$
\sn
\item  $\varepsilon^{\up}_\delta =
  \sup\{F'_1  
  (\gamma_{\ell}  
    (\zeta,\delta))(\ell)+1:\zeta \in t_\delta$ and $\ell
< k(\zeta,\delta)\} = 
\cup  \{\sup \Rang(F_1(\rho_{h}
(\zeta,\delta))) +1:\zeta \in t_\delta\}$; 
as $\cf(\delta) = \partial $ 
and $\extheta  = \cf(\extheta ) > |t_\delta|$, necessarily
$\alpha_{1,\delta} < \delta$ and $\varepsilon^{\up}_\delta < \extheta $.
\end{itemize}
\mn
Lastly, there are $\alpha^*_1 < \lambda$ and
$\varepsilon^{\up} < \kappa_1 = \extheta $ and 
$\cU^{\up}_1 \subseteq S^*_0$ as required by using Fodor lemma.  So
$(*)_3$ holds indeed.]

\bigskip    

Now let $E = \{\delta < \lambda:\delta$ is a limit ordinal $>
\alpha^*_1$ such that $\delta = \sup(\cU^{\dn}_1 \cap \delta)$
and $\alpha < \delta \Rightarrow t_\alpha \subseteq \delta\}$, it is a
club of $\lambda$ because $\alpha^*_1 < \lambda$ by $(*)_3(c)$ and
$\cU^{\dn}_1$ is an unbounded subset of $\lambda$ by $(*)_2(a)$, 
and $ t_ \alpha $ is a subset of $ \lambda $ of cardinality 
$ < \extheta  $  hence is bounded.  

Choose $\varp(*) = \max\{\varp^{\up} +1,\varp^{\dn} +1,\varp + 1\}$
where $\varp$ is from $\boxplus(c)$,
so $\varp(*) < \extheta $ and choose $\delta_2 \in E \cap S$ such that
$F'_1(\delta_2) = \varp(*)$.  Next choose $\alpha_2 \in \cU^{\up}_1
\backslash (\delta_2 +1)$ and let $\alpha^* \in (\alpha^*_1,\delta_2)$ 
be large enough such that $\zeta \in
(\alpha^*,\delta_2) \wedge \xi \in t_{\alpha_2} \Rightarrow
\rho(\xi,\delta_2) \char 94 \langle \delta_2 \rangle \triangleleft
\rho(\xi,\zeta)$.  Now choose $\delta_1 \in \cU^{\dn}_1 \cap
(\alpha^*,\delta_2)$ and $\alpha^{**} \in (\alpha^*,\delta_1)$ be such
that $\alpha \in (\alpha^{**},\delta_1) \wedge \xi \in t_{\alpha_2}
\Rightarrow \rho(\xi,\delta_1) \char 94 \langle \delta_1\rangle 
\triangleleft \rho(\xi,\alpha)$.

[Why this is possible?   
First as $ \alpha ^{**} > \alpha ^*$   it is enough to have
$ \alpha \in (\alpha ^{**}, \delta _1)  
    \Rightarrow 
    \rho ( \delta _2 , \delta _1 ) {\char 94} \langle \delta _1\rangle 
    \triangleleft \rho ( \delta _2, \alpha ) $. 
    Second here $ \cf(\delta _1) < \partial $ 
    however  this condition holds  because 
    $ \delta _ 1 \in {\mathscr U } ^{dn}_1 \subseteq \W$ so necessarily
    $ \gamma _{\lt}(\delta _2, \delta _1)= \delta _1 + 1 $ by 
    $ \odot _5(vi)$].  

Next let $\ell_* < \ell g(\rho(\alpha_2,\delta_1)$ be such that:
\mn
\begin{enumerate} 
\item[$(*)_4$] 
\begin{itemize}
\item[(a)] 
$  \varepsilon ( \bullet )   
    := F_1(\rho_{h}(\alpha_2,\delta_1))(\ell_*) = \max \Rang
  F_1(\rho_{h}(\alpha_2,\delta_1))$
\sn
\item[(b)]   under this restriction $\ell_*$ is minimal.
\end{itemize}
\end{enumerate}  
\mn

Lastly, choose $\alpha_1 \in (\alpha^{**},\delta_1)$ which is as in
$(*)_2(d)$ with respect to $\delta_1$, i.e. such that:
\mn
\begin{enumerate}
\item[$(*)_5$]  if $\zeta \in t_{\alpha_1}$ then $\Rang F_1(\rho_{h}
(\delta_1,\zeta)) \subseteq \varp^{\dn}$.
\end{enumerate}
\mn
Now we shall prove that the pair $(\alpha_1,\alpha_2)$ is as
required.  So let $(\zeta,\xi) \in t_{\alpha_1} \times t_{\alpha_2}$;
now by our choices
\mn

\begin{enumerate}
\item[$(*)_6$]  $\rho(\xi,\zeta) = \rho(\xi,\alpha_2) \char 94 
\rho(\alpha_2,\delta_2) \char 94 \rho(\delta_2,\delta_1) \char 94 
\rho(\delta_1,\zeta)$  
and $ \rho ( \alpha _2, \delta _1) =
\rho(\alpha_2,\delta_2) \char 94 \rho(\delta_2,\delta_1) $   
\end{enumerate}
\mn
So
\mn
\begin{enumerate}
\item[$(*)_7$]  $\Rang(F_1(\rho_{h}(\xi,\alpha_2)) \subseteq
 \varp^{\up} \le \varp(*)$
 \end{enumerate} 
\sn
[Why? by $ (*)_3$(a), the choice of $ \alpha _2 \in {\mathscr U } ^{\up}_1$
and $ \xi $ being from $ t_{\alpha _2}$]
\begin{enumerate} 
\item[$(*)_8$]  $\Rang(F_1(\rho_{h}(\delta_1,\zeta)) \subseteq
 \varp^{\dn} \le \varp(*)$
\end{enumerate} 
\sn
[Why by $ (*)_2$(d) and the choice of $ \alpha _1$ (and
$ \zeta $ being a member of $ t_{\alpha _1}$]
\begin{enumerate} 
\item[$(*)_9$]  $\varp(*) = F_1 \circ h(\delta_2) \in
  \Rang(F_1(\rho_{h}(\alpha_2,\delta_1)))$, see 
  $ (*)_6$  and 
  (before and
  after) $\odot_1$ .
\end{enumerate}
\mn
[Why?  Recall that $\delta_2$ was chosen in $E \cap S$ such that
$F'_1(\delta_2) = \varp(*)$.]

Hence 

\begin{enumerate}    
\item[$(*)_{10}$] $\varepsilon 
            \le \varepsilon (*)  \le \varepsilon ( \bullet ) < \partial     $
\end{enumerate} 

Putting those together, We can finish this stage by:   
\mn
\begin{enumerate}
\item[$(*)_{11}$]  in $\boxplus_3(b)$ for our $\bar t$ and the pair
  $(\alpha_1,\alpha_2)$, our $  \varepsilon ( \bullet ) $ 
  (chosen before $(*)_5$) is
  gotten, witnessing $  \varepsilon ( \bullet )  
  \in A_{\bar t,\varp(*)} \subseteq
  A_{\bar t,\varp}$ 
  \end{enumerate} 
  [Why?  
  As first  
  $\varp < \varp(*)$, 
  by the choice of $ \varepsilon (*)$,
  and second 
if $(\zeta,\xi) \in t_{\alpha_1} \times t_{\alpha_2}$ then $\ell = 
\ell g(\rho(\xi,\alpha_2)) + \ell_*$ is as required in $\boxplus_3(b)$
for $\bar t  $   
by $ (*)_6 - (*)_{10}$]  
\mn
So we are done proving $\boxplus_7(a)$.  
\medskip

\noindent
\underline{Stage D}:  
By $\boxplus_8$
\mn
\begin{enumerate}
\item[$\circledast_1$]   there is $F_*:\lambda \rightarrow \extheta $ such that
$\varepsilon < \extheta \Rightarrow F^{-1}_*(\{\varepsilon\}) \ne
  \emptyset \mod D$.
\end{enumerate}
\mn
We first deal with the easier version with
$\extheta $ colours, i.e. proving $\Pr_1(\lambda,\lambda,\extheta ,\extheta )$.

We now define the colouring $\bfc_1:[\lambda]^2 \rightarrow \extheta $ by:
\mn
\begin{enumerate}
\item[$\circledast_2$]  if $\alpha < \beta < \lambda$ then 
$\bfc_1\{\alpha, \beta \}$ is 
$F_*(\gamma_{\ell(\beta,\alpha)}(\beta,\alpha))$  
  where $\ell(\beta,\alpha) = \min\{\ell <
  k(\beta,\alpha):F'_1(\gamma_\ell(\beta,\alpha)) = \max
\Rang(F'_1(\rho(\beta,\alpha)))\}$.
\end{enumerate}
\mn
To prove that the colouring $\bfc_1$ really witnesses
$\Pr_1(\lambda,\lambda,\extheta ,\extheta )$, our task is to prove:
\mn
\begin{enumerate}
\item[$\circledast_3$]  given $\bar t \in \bfT$ and $\iota < \extheta $ there are
  $\alpha <  \beta$ such that:
\sn
\begin{itemize}
\item  $\zeta \in t_\alpha \wedge \xi \in t_\beta \Rightarrow
  \bfc_1\{\zeta,\xi\} = \iota$.
\end{itemize}
\end{enumerate}
\mn
[Why does $\circledast_3$ holds? 
Let $B_\iota = \{\gamma < \lambda:F_*(\gamma) = \iota\}$.  By
the choice of $F_*$ we know that $B_\varp \ne \emptyset \mod D$.
Focus on $A_{\bar t,\varp}$ for our   
specific $\bar t \in \bfT$ and any
$\varp < \extheta $.  Since 
$A_{\bar t,\varp} \in D$ we conclude that $B_\varp \cap A_{\bar t,\varp} \ne
\emptyset$.

Fix an ordinal $\gamma \in B_\iota \cap A_{\bar t,\varp}$.  By the
very definition of $A_{\bar t,\varp}$ in $\boxplus_3$ we choose
$\alpha < \beta < \lambda$ 
such that for every $(\zeta,\xi) \in
t_\alpha \times t_\beta$ there exists $\ell < k(\xi,\zeta)$ for which
$\gamma_\ell(\xi,\zeta) = \gamma$ and $F'_1(\gamma) \ge
F'_1(\gamma_k(\xi,\zeta))$ whenever $k < k(\xi,\zeta)$ and
$F_1(\gamma) \ge \varp$ and $F'_1(\gamma) > F'_1(\gamma_k(\xi,\zeta))$
whenever $k < \ell$.  Let $\ell(\xi,\zeta)$ be this
$\ell$, in fact, this $\ell$ is unique (for the pair $(\zeta,\xi)$).

Now $\bfc_1\{\zeta,\xi\} = F_*(\gamma_{\ell(\xi,\zeta)}(\xi,\zeta))$ (by
$\circledast_2$) which equals $F_*(\gamma)$ (by the choice of
$\ell(\xi,\zeta)$) which equals $\iota$ (since $\gamma \in B_\iota$).
Hence $\circledast_3$ holds and we finish Stage D.]
\medskip

\noindent
\underline{Stage E}:  The full theorem: the case of $\lambda$ colors.

Let $h',h''$ be functions from $\extheta $ into $\extheta ,\omega$
respectively such that the mapping $\zeta \mapsto 
(h'(\zeta),h''(\zeta))$ is onto $\extheta \times \omega$ and moreover
each such pair is gotten $\extheta $ times.

We have to define a colouring $\bfc_2:[\lambda]^2 \rightarrow \lambda$
exemplifying $\Pr_1(\lambda,\lambda,\lambda,\extheta )$.

This is done as follows using $h',h''$ and $F_*$ from $\circledast_1$:
\mn
\begin{enumerate}
\item[$\oplus_1$]  for $\alpha < \beta < \lambda$ we let
\sn
\begin{enumerate}
\item[$\bullet_1$]   $\zeta = \zeta(\beta,\alpha) :=
  h'(\bfc_1\{\beta,\alpha\})$, necessarily $<\extheta $
\sn
\item[$\bullet_2$]  $n = n(\beta,\alpha) :=
  h''(\bfc_1\{\beta,\alpha\})$, necessarily $< \omega$
\sn
\item[$\bullet_3$]  $m = m(\beta,\alpha)$ is the $n$-th member of $\{k
  < k(\beta,\alpha):F'_1(\gamma_k(\beta,\alpha))=\zeta\}$ when there is
  such $m$ and is zero otherwise
\sn
\item[$\bullet_4$]  we define $\bfc_2$ as follows: for $\alpha <
  \beta,\bfc_2\{\alpha,\beta\}$ is
  $F'_2(\gamma_{m(\beta,\alpha)}(\beta,\alpha))$ 
     recalling that $ F'_2 $, a function from $ \lambda$ 
     to $\lambda $ 
     is from $ \odot_2$  from
     the end of stage A.
\end{enumerate}
\end{enumerate}
\mn
To prove that $\bfc_2$ indeed exemplifies
$\Pr_1(\lambda,\lambda,\lambda,\extheta )$   it suffice 
to prove (this is the
task of the rest of the proof)
\mn
\begin{enumerate}
\item[$\oplus_2$]  assume $\bar t \in \bfT$ and $j_* < \lambda$ and we
  shall 
  find $\alpha < \beta$ such that $t_\alpha \subseteq \beta$
and $(\zeta,\xi) \in t_\alpha \times t_\beta \Rightarrow
  \bfc_2\{\zeta, \xi \} = j_*$.  
\end{enumerate}
\mn
Toward this:
\mn
\begin{enumerate}
\item[$\oplus_3$]
\begin{enumerate}
\item[(a)]  we apply $(*)_3$ to our $\bar t$, getting
  $\varp^{\up},\cU^{\up}_1,\alpha^*_1$ as there
\sn
\item[(b)]  we apply $(*)_2$ to our $\bar t$ getting
  $\cU^{\dn}_1,\varp^{\dn}$ 
\sn
\item[(c)]  let $\varp^{\md} = \max\{\varp^{\up} +1,
  \varp^{\dn} +1\}$.
\end{enumerate}
\end{enumerate}
\mn
We can find $g_2,\cU^{\up}_2,\gamma _*,\alpha^*_2,m^*_2$ 
such that:
\mn
\begin{enumerate}
\item[$\oplus_4$]
\begin{enumerate}
\item[(a)]  $\gamma_* < \lambda$ satisfies  
  $F_2(\gamma_*) = j_*$ and $F_1(\gamma_*) = \varp^{\md}$
\sn
\item[(b)]  $\cU^{\up}_2 \subseteq S^*_{\gamma _*}$ is stationary 
hence  
  $\delta \in \cU^{\up}_2 \Rightarrow F'_2(\delta) = F_2(h(\delta)) =
  F_2(\gamma_*) = j_* \wedge F'_1(\delta) = F_1(h(\delta)) =
  F_1(\gamma_*) = \varp^{\md}$
\sn
\item[(c)]   $g_2$ is a function with domain 
$\cU^{\up}_2$ such that $\delta \in \cU^{\up}_2
  \Rightarrow \delta < g_2(\delta) \in \cU^{\up}_1$
\sn
\item[(d)]  $\alpha^*_2$ satisfies 
$\alpha^*_1 < \alpha^*_2 < \min(\cU^{\up}_2)$ 
\sn
\item[(e)]  if $\delta \in \cU^{\up}_2$ and $\alpha \in
[\alpha^*_2,\delta)$ and $\beta \in t_{g_2(\delta)}$ then
\sn
\begin{itemize}
\item   $\rho(g_2(\delta),\delta)
\char 94 \langle \delta \rangle \trianglelefteq \rho(g_2(\delta),\alpha)$ hence
\sn
\item  $\rho_{\beta,\delta} \char 94 
\langle \delta \rangle \trianglelefteq \rho_{\beta,\alpha}$
\end{itemize}
\sn
\item[(f)]   $m^*_2$ satisfies:  
for every $\delta \in \cU^{\up}_2$,  
it is  
the cardinality of the set $\{\ell < k(g_2(\delta),\delta):
F'_1(\gamma_\ell(g_2(\delta),\delta)) = \varp^{\md}\}$ which may be zero.
\end{enumerate}
\end{enumerate}
\mn
[Why?  First choose $\gamma _*$ as in clause (a) of $\oplus_4$
(possible by the choice of 
 $ F_1,F_2$ in the beginning of Stage A;
hence 
$ \delta \in S_{\gamma_* } \Rightarrow F'_2(\delta )=
F_2(h(\delta ))   =  
F_2 ( \gamma _*)  
= j_*$   
and $ F'_1( \delta )= F_1(h( \delta ))= F_1(\gamma _*)= \varp^{\md}$ 
(by the choice of 
$F'_1$ in $\odot_7$ recalling the definitions of $h,F'_1$).  
Second, define $g':S^*_{\gamma _*} \rightarrow \cU^{\up}_1$ 
such that $\delta \in S^*_{\gamma _*} \Rightarrow \delta <
 g'(\delta) \in \cU^{\up}_1$.  
Third, for each $\delta \in S^*_{\gamma _*} \backslash (\alpha^*_1+1)$, 
find $\alpha'_{2,\delta} < \delta$ above $\alpha^*_1$
and $m_{2,\delta}$ such that the parallel of clauses 
(e),(f)
(with $ g' $ here instead of $ g_2 $ there) 
of $\oplus_4$ holds.  Fourth, use Fodor lemma to get a
stationary $\cU^{\up}_2 \subseteq S^*_{\gamma _*}$ such that $\langle
(\alpha'_{2,\delta},m_{2,\delta}):\delta \in \cU^{\up}_2\rangle$ is
constantly $(\alpha^*_2,m^*_2)$ and lastly let
$g_2 = g' \rest \cU^{\up}_2 \backslash (\alpha^*_2 +1)$.   
   Now it is easy to check  that   
$\oplus_4$ holds indeed.]

Next
\mn
\begin{enumerate}
\item[$\oplus_5$]  if $\delta \in \cU^{\up}_2$ \then \,:
\sn  
\begin{enumerate}
\item[(a)]  $F'_1(\delta) = \varp^{\md}$
\sn
\item[(b)]  if $\alpha \in [\alpha^*_2,\delta),\xi \in t_{g_2(\delta)}$
  \then \, $u = \{\ell < k(\xi,\alpha):F'_1(\gamma_\ell(\xi,\alpha)) =
  \varp^{\md}\}$ has $> m^*_2$ members and if $\ell$ is the $m^*_2$-th
  member of $u$ then $\gamma_\ell(\xi,\alpha) = \delta$.
\end{enumerate}
\end{enumerate}
\mn
Why?  Clause (a) holds by $\oplus_4(a),(b)$.  For clause (b) use
clause (a) and the demands on $m^*_2$.
That is   

\begin{enumerate}
\item[(a)] $ \rho (  \xi ,  \alpha  ) = \rho ( \xi , g_2( \delta ) )
  {\char 94} \rho ( g_2( \delta ) , \delta ) 
  {\char 94} \rho (\delta, \alpha ) $  
 
 [Why? by $ (*)_3, \oplus _4(e)$]
\item[(b)] $ \Rang (\rho _{h}(\alpha, g_2 (\delta ))) \subseteq 
   \varepsilon ^{\up}  \subseteq \varepsilon ^{\md}$
   
   [Why? by $ (*)_2$]
\item[(c)] the set
   $ \{ {\ell} < k(g_2(\delta ) , \delta 
   ):  F'_1( \gamma _ {\ell} 
       (g_2(\delta), \delta )) = \varepsilon ^{\md } \} $ 
       has $ m^*_2 $ members 
       
       [why?  by $ \oplus _4(f)$]
\item[(d)] $ F'_1 ( \gamma _0 ( \delta, \alpha ) )= 
       F'_1(\delta ) = \varepsilon ^{\md }$
       
       [Why? by $ \oplus _4(a),(b)$]]
\item[(e)] 
if $ {\ell} _* $ is the $ m^*_2 $-th member 
of $ \{ {\ell} : F_1( \gamma _ {\ell}   
(\xi , \alpha ) )= \varepsilon ^{\md }\} $
    then $ \gamma _{{\ell} _*}(\xi, \alpha )= \delta $ 
    
    [Why? putting the above together]   
  \end{enumerate} 

So $ \oplus _5$  holds indeed. 


Now choose $\varp(*) < \extheta $ such that $h'(\varp(*)) = \varp^{\md}$
and $h''(\varp(*)) = m^*_2$.

Next, let $E = \{\delta < \lambda:\delta$ limit ordinal $> \alpha^*_2$
such that $\delta = \sup(\cU^{\dn }_1
\cap \delta) $ 
and $\alpha < \delta
\Rightarrow  
g_2(\alpha) < \delta\}$. 

Lastly, 

\begin{enumerate} 
\item[$ \oplus_6$]
choose $\delta_1 < \delta_2$ such that 
\begin{enumerate}
\item[(a)] $\delta_1 \in
\cU^{\dn}_1 \cap E$ 
\item[(b)] $\delta_2 \in \cU^{\up}_2 \cap E \backslash
(\delta_1 + 1)  $ 
\item[(c)] $ \bfc_1\{\delta_2,\delta_1 \}  = \varp(*)$, 
\end{enumerate} 
\end{enumerate} 

[Why does such a pair $ (\delta _1, \delta _2  ) $ exist? 
By Stage
D applied to $\bar s = \langle s_\alpha:\alpha < \lambda\rangle$ where
$s_\alpha = \{\min(\cU^{\dn}_1   
    \cap E \backslash
\alpha),\min(\cU^{\up}_2 \cap E \backslash \alpha)\}$.

That is, we can find ordinals $ \alpha < \beta  < \lambda $
such that:
for every $ (\zeta, \xi ) \in (s_{\alpha } \times s_\beta )$
we have $ \mathbf{c} _1 \{  \xi, \zeta \} = \varepsilon ^{\md }$.

Let 
$ \delta_1 = \min ( {\mathscr U } ^{\dn}_1 \cap E \setminus \alpha $
and let 
$ \delta_2 = \min ( {\mathscr U } ^{\up}_1 \cap E \setminus \beta  $.

So $ (\delta _1, \delta  _2 ) \in (s _ \alpha \times s_\beta )$
hence clearly $ \delta _1 < \delta _2$, 
$ \mathbf{c} _1 \{ \delta _1 , \delta _2  \}  = \varepsilon (*)$,
$ \delta _1  \in {\mathscr U } ^{\dn}_1  \cap E$
and $ \delta _  1   
    \in {\mathscr U } ^{\up}_1 \cap E $.
So the pair $ (\delta _1, \delta _2 )$ is as promised in
in $  \oplus _6$]


Now let $\beta = g_2(  
\delta_2)$ and choose $\alpha
\in \cU^{\dn}_1   
\cap \delta_1 \backslash (\alpha^*_2 +1)$.  Easy to
check that $\alpha,\beta$ are as required.

So we have finished proving Theorem \ref{f4}.
\end{PROOF}

\bibliographystyle{amsalpha}
\bibliography{shlhetal}

\end{document}